\documentclass[conference,letterpaper]{IEEEtran}

\usepackage[utf8]{inputenc} 
\usepackage[T1]{fontenc}
\usepackage{url}
\usepackage{ifthen}
\usepackage{cite}
\usepackage[cmex10]{amsmath} 

\def\withnotes{1}
\def\withcolors{1}
\usepackage{tex_macros}
\usepackage{graphicx}
\usepackage{subcaption}
\usepackage{bm}
\usepackage{framed}
\usepackage{physics}
\usepackage{nicefrac}
\usepackage[most]{tcolorbox}
\usepackage{mathrsfs}
\tcbset{colback=blue!10!white, colframe=blue!50!black, 
        highlight math style= {enhanced, %<-- needed for the ’remember’ options
            colframe=blue,colback=red!10!white,boxsep=0pt}
        }

\def\Zpart{\mathcal{Z}}
\def\R{\mathbb{R}}
\def\E{\mathbb{E}}

\def\one{\bm{1}}

\def\Gauss{{\cal N}}

\def\fe{{\cal F}}

\begin{document}
\title{Analysis of Diffusion Models for Manifold Data} 

\author{%
  \IEEEauthorblockN{Anand Jerry George, Rodrigo Veiga, Nicolas Macris}
  \IEEEauthorblockA{EPFL, School of Computer and Communication Sciences.\\
                    CH-1015 Lausanne, Switzerland.
                    }

}

\maketitle

\begin{abstract}
   We analyze the time reversed dynamics of generative diffusion models. If the exact empirical score function is used in a regime of large dimension and exponentially large number of samples, these models are known  to undergo transitions between distinct dynamical regimes. We extend this analysis and compute the transitions for an analytically tractable manifold model where the statistical model for the data is a mixture of lower dimensional Gaussians embedded in higher dimensional space. We compute the so-called speciation and collapse transition times, as a function of the ratio of manifold-to-ambient space dimensions, and other characteristics of the data model. An important tool used in our analysis is the exact formula for the mutual information (or free energy) of Generalized Linear Models.
\end{abstract}

\section{Introduction}

In generative modeling, we are concerned with the following problem. Given a set $S$ of i.i.d. samples $\{x_i\}_{i=1}^n$ in $\mathbb{R}^d$,  from an {\it unknown} probability distribution $\pi$, we want to generate a new sample from $\pi$ independent of $S$. Generative diffusion models \cite{sohl-dickstein_deep_2015,song_generative_2019,ho_denoising_2020,yang_diffusion_2024} have emerged as an interesting tool for this task. These models leverage stochastic processes guided by a \emph{score function} to iteratively transform simple distributions, such as Gaussian noise, into non-trivial data distributions such as $\pi$ \cite{anderson_reverse-time_1982, haussmann_time_1986}. In practice, the optimal score function is unknown (because $\pi$ is unknown) and has to be estimated from the sample set $S$. However, it is theoretically unclear how to achieve this so that good generalization is achieved, as opposed to mere memorization~\cite{yoon_diffusion_2023}.
More generally, diffusion models seem to undergo distinct dynamical transitions in their behaviors \cite{biroli_generative_2023,biroli_dynamical_2024,raya_spontaneous_2023,li_critical_2024}, whose comprehensive understanding remains incomplete. Given this state of affairs, it is of theoretical value to explore the dynamical properties of diffusion models using a naive empirical score function. 

Considering data generated by a mixture of Gaussians in $\R^d$, Ref.~\cite{biroli_dynamical_2024} identified three distinct dynamical behaviors of the backward generative diffusion process with empirical score, in the regime $n=e^{\alpha d}$, $d\to \infty$, and $\alpha$ fixed. First, the reversed process starts from {\it pure noise} and the random trajectories do not capture any data structure. Second, after a {\it speciation time}, the reverse trajectory \emph{specializes} into one of the two classes of the data. Third, after a {\it collapse time}, trajectories are confined to the basins of attraction of a data points and \emph{collapse} towards them. Speciation occurs on a time scale $t_S\sim \log d$ and collapse time $t_C= \frac{1}{2}\log(1+ (e^{2\alpha}-1)^{-1})$ corresponds to a sharp transition in the limit $d\to +\infty$, $(\nicefrac{1}{d})\log n \to \alpha$. 

It is of interest to investigate settings with structured data reflecting real datasets such as images, text, etc, using diffusion models~\cite{sclocchi_phase_2025}. 

In this work, we focus on a simple tractable model of structured data. We consider a statistical model representing data as a mixture of lower $p$-dimensional Gaussians in a manifold embedded in a higher $d$-dimensional ambient space ($d>p$). This is motivated by the observation that real-world high-dimensional data often effectively resides on lower-dimensional manifolds. Here, for the manifold, we take a $p$-dimensional hyperplane which is then warped by applying a point-wise non-linear function (e.g., a sigmoid activation). Such manifold models have already been used in the learning theory and inference context where they provide a tractable setting (see, e.g.,~\cite{goldt_modeling_2020,hand_global_2020,luneau_tensor_2020}). Closer to this work, Refs.~\cite{achilli_losing_2024,ventura_manifolds_2024} have investigated the dynamical regimes in diffusion models, when the data lies on a linear manifold. We remark that for the case of linear manifolds, our result on collapse time is consistent with the results of~\cite{achilli_losing_2024}.

Using the empirical score, we derive in Sections \ref{sec:Stime} and \ref{sec:Ctime} explicit results for the speciation and collapse times $t_S, t_C$, in the regime $d, p, n\to +\infty$ with $p=\beta d$ and $n= e^{\alpha d}$ for fixed $\alpha>0,\; 0<\beta<1$. The setting is introduced in Section~\ref{setting} and our main contributions are summarized in Section \ref{subsec:contributions}.

\section{Diffusion models for data in a manifold}\label{setting}

Diffusion models solve the generative modeling problem by time reversing a diffusion process that transports $\pi$ to a known distribution such as an isotropic Gaussian~\cite{sohl-dickstein_deep_2015,song_generative_2019,ho_denoising_2020}. Consider the forward $d$-dimensional Ornstein-Uhlenbeck process \cite{gardiner_stochastic_2009} (with standardized variance)
$\dd X_t = - X_t\dd t +\sqrt{2}\;\dd W_t$ with $X_0\sim\pi$.
The conditional distribution of $X_t$ given $X_0$ is given by a Gaussian distribution $\Gauss(a_t X_0,h_tI_d)$, where $a_t=e^{-t}$ and $h_t=1-e^{-2t}$. The probability distribution of $X_t$ is
\begin{equation}
    P_t(x) = (2\pi h_t)^{-d/2}\int_{\R^d}  \dd x_0 \; e^{-\frac{\norm{x-a_tx_0}^2}{2h_t}}\pi(x_0)\;,
\end{equation}
though $\pi$ is unknown. The time reversed process satisfies the following stochastic differential equation\cite{anderson_reverse-time_1982}:
\begin{align}\label{backward}
    -\dd Y_t = (Y_t+2\nabla\log P_{t}(Y_t)) \;  \dd t+\sqrt{2} \; \dd W_t \;,
\end{align}
which runs backward in time starting from $Y_T\sim P_T$. Here $P_T$ is unknown, but for $T$ large is very close to $\cN{0,I_d}$, with $I_d$ denoting the $d\times d$ identity matrix. So, we start the reverse process with $Y_\infty\sim\cN{0,I_d}$ without incurring much error. It is a well-known old result that the backward process converges to $Y_0\sim \pi$ \cite{haussmann_time_1986}: if the so-called {\it score function} $s(x, t)=\nabla\log P_{t}(Y_t)$ were known, we could use the dynamics to sample from $\pi$.

The learning task is to estimate $s(x, t)$ using the set of samples $S=\{x_1,x_2.\cdots,x_n\}$. The most naive choice is to estimate $\pi(x)$ by the empirical distribution $\frac{1}{n}\sum_{i=1}^n\delta(x- x_i)$ and take the {\it empirical score}: $s^e(t,x) = \nabla\log P_t^e(x)$, where $P_t^e(x)$ is given by
\begin{equation}
    \label{eq:pte}
    P_t^e(x) = n^{-1}(2\pi h_t)^{-d/2}\sum_{i=1}^ne^{-\frac{\norm{x-a_tx_i}^2}{2h_t}} \;.
\end{equation}
As shown in \cite{song_maximum_2021}, this is also the minimizer of an appropriate quadratic loss function. In this paper we are concerned with the dynamical regimes induced by this empirical score function. 

Our model for the data samples is as follows. The samples in ambient space $x_i\in \R^d$, $i=1,\dots, n$ are assumed to lie in a lower dimensional manifold $x_i = \phi(\frac{F\xi_i}{\sqrt p})$, $\xi_i\in \mathbb{R}^d$, $p< d$, $F$ is a $d\times p$ matrix with real entries, and $\phi$ an activation function acting component-wise. The matrix $F$ will be taken with random i.i.d $\mathcal{N}(0, 1)$ entries or with a set of $p$ orthogonal columns. The data points in the lower dimensional manifold are sampled from a simple mixture of two Gaussians with p.d.f. $q(\xi) = \frac{1}{2} q_+(\xi) + \frac{1}{2} q_-(\xi)$, where  
\begin{equation}
    \label{eq:qpm}
    q_\pm (\xi) = (2\pi\rho)^{-p/2}e^{-\frac{\Vert\xi - \mu_{\pm}\Vert^2}{2\rho}} \;,
\end{equation}
for $\rho >0$ and $\mu_{\pm}\in\R^p$.

When the activation is linear $\phi(u) = u$, the data lie in a $p$-dimensional hyperplane. If, furthermore, $d=p$ and $F/\sqrt p$ is the identity matrix, the basic model studied in \cite{biroli_dynamical_2024} is recovered.

\subsection{Summary of main contributions}\label{subsec:contributions}

We look at a regime of large dimensions and exponentially large number of samples. More precisely $d, p\to +\infty$, $p/d= \beta$, $n= e^{\alpha d}$, $\alpha>0$ and $0<\beta<1$ fixed.

In Section \ref{sec:Stime} we analyze the specialization phenomenon. In this short note, we carry out the details for the simplest case of opposite centers $\mu_+ + \mu_-=0$ and odd activation functions. More general cases may be approached by the same methods but would require much more elaborate analysis and discussion. We find that the effect of the non-linearity is entirely captured by the quantity $\Gamma_0 (y)\equiv \E_{u\sim\Gauss(0,1)}\left[ \phi\left(\sqrt{\rho}\;u+y\right)\right]$. Let
 $\varrho_1 =\E_{u\sim\Gauss(0,1)}\left[\Gamma_0 (u)u\right]$, $
 \varrho_*^2 =\E_{u\sim\Gauss(0,1)}\left[\Gamma_0 (u)^2\right] - \varrho_1^2$.  Let also 
$\tilde{\mu}_\pm=\mu_\pm/\sqrt{p}$ be the normalized center of the mixtures.
We find the expression   
$t_S \approx \frac{1}{2} \log[ 2 ( \varrho_1^2  \beta d  \norm{\tilde{\mu}_\pm}^2  + \varrho_*^2  )]
$ (valid for $p$ and $n$ large).
For the case of a linear manifold $\varrho_1=1$, $\varrho_*=0$ and the formula simply reduces to $\frac{1}{2}\text{log}(\beta d \norm{\tilde{\mu}_\pm}^2)$. For $p=d$ we recover the expression of \cite{biroli_dynamical_2024}.

In Section~\ref{sec:Ctime} we investigate the collapsing regime. We follow the approach of~\cite{biroli_dynamical_2024} using an analogy with the Random Energy Model (REM). For times $t<t_C$ (the end of the reversed process corresponding to a collapsing phase), the empirical distribution \eqref{eq:pte} along a trajectory of the process is dominated by one data sample, say the term $i=1$. For $t>t_C$ on the other hand, it is the rest of the sum for $i\geq 2$ which dominates, and is well approximated by the partition function of a Generalized Linear Model (GLM). Using the exact formula for the free energy (average of log-partition function or mutual information) of the GLM \cite{barbier_optimal_2019}, we can compute $t_C$ through Eq.~\eqref{eq:TC}, which involves only one-dimensional integrals and optimization of a function involving two scalar parameters. For a linear activation, Eq.~\eqref{eq:TC} reduces to $t_C= \frac{1}{2}\log(1+ (e^{2\alpha/\beta}-1)^{-1})$ and for $p=d$ we get back the result in \cite{biroli_dynamical_2024}. As shown in this reference, $t_C$ also corresponds to the condensation phase transition of the REM, and this extends to the present manifold model. Thus, in the asymptotic limit of infinite dimension, the change of dynamical behavior of the reversed process is a sharply defined transition at $t_C$. 

\section{Speciation time}\label{sec:Stime}

The speciation transition occurs at the beginning of the backward dynamics, for large times. In this regime, $h_t$ is exponentially close to $1$ and $a_t$ is exponential small. Therefore, for large $t$ (and fixed separation between the centers $\mu_\pm$) the distributions $P_t^e$ and $P_t$ are not very different. In this regime, we replace the empirical score function $s^e(x,t)$ by the exact score function $s(x, t)$. For the mixture of two Gaussians, the exact score can be written as
\begin{equation}
\label{eq:scoregm}
\nabla\log P_{t}(x) = \frac{1}{2} \nabla \log \left( z_t^{+}(x) + z_t^{+}(x) \right) \;,
\end{equation}
with $z_t^\pm(x) = e^{- \frac{1}{2h_t}\beta x^\top x+g_{\pm}(x)}$ and
\begin{equation}
    g_{\pm}(x) = \log\E_{\xi_\pm} \left[ e^{\frac{a_t}{h_t} \phi\left(\frac{F\xi}{\sqrt p}\right)^\top x -\frac{a_t^2}{2h_t} \phi\left(\frac{F\xi}{\sqrt p}\right)^\top \phi\left(\frac{F\xi}{\sqrt p}\right)} \right] \;,
\end{equation}
where $\E_{\xi_\pm}$ indicates expectation over $q_{\pm}(\xi)$ given in Eq.~\eqref{eq:qpm}. In the limit of large times, we can expand this expression around $a_t$~\cite{biroli_dynamical_2024}:
\begin{align}
    g_{\pm} (x) &=\frac{a_t}{h_t} \sum_{j=1}^d x_j \zeta_j^\pm  + \frac{a_t^2}{2 h_t} \sum_{j=1}^d \left[ \left(x_j^2 - h_t \right) \zeta_{jj}^\pm - x_j^2 (\zeta_j^\pm)^2  \right] \nonumber \\ 
      &+  \frac{a_t^2}{h_t} \sum_{j=1}^d \sum_{l\neq j}^d  \left[ \zeta_{jl}^\pm - \zeta_j^\pm \zeta_l^\pm  \right]  \;,
\end{align}
where
\begin{align}
    \zeta_j^\pm    &\equiv \E_{z}\left[\phi\left( \sqrt{\rho} \left(\nicefrac{f_j^\top z}{\sqrt p}\right)+ \lambda_j^\pm \right)\right] \;,  \nonumber \\
    \zeta_{jl}^\pm &\equiv \E_{z}\left[\phi\left( \sqrt{\rho} \left(\nicefrac{f_j^\top z}{\sqrt p}\right) + \lambda_j^\pm \right)\phi\left( \sqrt{\rho} \left(\nicefrac{f_l^\top z}{\sqrt p}\right) + \lambda_j^\pm \right)\right]  \nonumber \;,
\end{align}
with $z\sim\Gauss(0,I_p)$ and $f_j^\top \in \R^p$ denoting the $j$-th row of $F$. The quantities $\lambda_j^\pm \equiv \frac{f_j^\top \mu_\pm}{\sqrt p}$ enclose the information about the centers of the Gaussian clouds. By the central limit theorem, when $p\to\infty$:
\begin{align}
     \zeta_{j}^\pm =&  \E_{u\sim\Gauss(0,1)}\left[\phi\left(\sqrt{\rho}\; u +\lambda_j^\pm\right)\right]  \;, \\
     \zeta_{jl}^\pm =&  \E_{u,v\sim\Gauss(0,\Theta_{jl})}\left[\phi\left(\sqrt{\rho}\; u +\lambda_j^\pm\right)\phi\left(\sqrt{\rho} v +\lambda_j^\pm\right)\right] \;,  
\end{align}
where $\Theta_{jl} \in  \R^{2\times 2}$ with matrix elements $\theta_{jl} = f_j^\top f_l / p $. In order to simplify the crossed-term for $j\ne l$, we perform an expansion in terms of Hermite polynomials using Mehler's formula~\cite{kibble_extension_1945}. Neglecting contributions of order $\nicefrac{1}{p}$:
\begin{equation}
    \zeta_{jl}^\pm = \Gamma_0\bigl(\lambda_j^{\pm}\bigl)\Gamma_0\bigl(\lambda_l^{\pm}\bigl)\;+\;\theta_{jl}\; \Gamma_1\bigl(\lambda_j^{\pm}\bigl)\Gamma_1\bigl(\lambda_l^{\pm}\bigl) \;,
\end{equation}
where
\begin{subequations}
\begin{align}
 \Gamma_0 (y)&\equiv \E_{u\sim\Gauss(0,1)}\left[ \phi\left(\sqrt{\rho}\;u+y\right)\right]    \;, \label{eq:gamma0} \\
 \Gamma_1 (y)&\equiv \E_{u\sim\Gauss(0,1)}\left[ \phi\left(\sqrt{\rho}\;u+y\right)u\right]    \;. \label{eq:gamma1} 
\end{align}
\end{subequations}
If $j=l$, this expansion is not useful, because the diagonal terms of $\Theta_{jl}$ tend to one and higher orders in $\theta_{jj}$ cannot be neglected. In any case, since we are interested in the dominant scaling of the speciation time, we write:
\begin{align}
    g_{\pm} (x) &=\frac{a_t}{h_t} \sum_{j=1}^d x_j  \Gamma_0\bigl(\lambda_j^{\pm}\bigl) +   \frac{a_t^2}{h_t} \sum_{j=1}^d \sum_{l\neq j}^d   \theta_{jl} \; \Gamma_1\bigl(\lambda_j^{\pm}\bigl) \Gamma_1\bigl(\lambda_l^{\pm}\bigl)  \nonumber \\
    & + \frac{a_t^2}{2 h_t} \sum_{j=1}^d \left[ \left(x_j^2 - h_t \right) \Gamma^{(2)}\bigl(\lambda_j^{\pm}\bigl) - x_j^2 \Gamma_1\bigl(\lambda_j^{\pm}\bigl)^2  \right]    \;,
\end{align}
with $\Gamma^{(2)}(y)\equiv\E_{u\sim\Gauss(0,1)}\left[ \phi\left(\sqrt{\rho}\;u+y\right)^2\right]$.

\subsection{Two equidistant Gaussians and odd activation function}
We consider the case of opposite centers and set $\mu_\pm = \pm \mu$ for fixed $\mu\in \R^p$. If the activation is an odd function, $\phi(y) = - \phi(-y)$, we have $\Gamma_0(\pm y)=\pm\Gamma_0 (y)$, $\Gamma_1 (\pm y) = \Gamma_1 (y)$ and $\Gamma^{(2)}(\pm y) = \Gamma^{(2)}(y)$. These symmetries imply cancellation of terms in the score and we find for the $j$-th component:
\begin{align}
    \label{eq:scorej}
    \partial_{x_j} \log P_t(x) =& -\frac{x_j}{h_t} + \frac{e^{-t}}{h_t}  \Gamma_0\bigl(\lambda_j\bigl) \tanh( e^{-t} \sum_{l=1}^d x_l  \Gamma_0\bigl(\lambda_l\bigl)  )  \nonumber \\
      & + e^{-2t} \Upsilon_{j} (x) \;,
\end{align}
where $\Upsilon_j (x)= (\nicefrac{x_j}{h_t^2})(\Gamma^{(2)}(\lambda_j) - \Gamma_0(\lambda_j)^2  ) + (\nicefrac{4}{h_t^2})  \sum_{l\ne j}^d  x_l  \; \theta_{jl} \; \Gamma_1 (\lambda_j) \Gamma_1 (\lambda_l)   $.

\begin{remark}
Note that a calculation for an even activation would show that the exponential factors proportional to $e^{-t}$ cancel and the leading order would be $e^{-2t}$. However, for opposite centers an even activation maps the centers $\pm \mu$ at the same point in ambient space and there is no speciation, so we do not discuss this case further. This remark becomes important for activations that have an even and odd part.
\end{remark}

Hereafter, we keep the leading contributions of order $e^{-t}$. Thus we neglect contributions proportional $e^{-2t}$, i.e., $\{ \Upsilon_j(x) \}_{j=1}^d$,  and also replace $h_t \approx 1$ to this same order. Within these approximations, by replacing Eq.~\eqref{eq:scorej} in the SDE~\eqref{backward}, we deduce that the scalar quantity $q \equiv \sum_{j=1}^d x_j \Gamma_0\left(\lambda_j\right)$, satisfies the stochastic equation:
\begin{equation}
    \label{eq:sdenlin}
    - \dd q = \left[ - q  + 2 e^{-t} \;   \sum_{\alpha=1}^d \Gamma_0\left(\lambda_\alpha\right)^2   \tanh( e^{-t}  q ) \right]\dd t + \dd \tilde{w}  \;, 
\end{equation}
where $ \dd \tilde{w}$ is the increment of a properly rescaled Wiener process.
Interpreting the drift term as a deterministic force given by the derivative of a \emph{potential}, this equation is $  - \dd q = - \pdv{q} V(q,t) \dd t + \dd\tilde{w}$, with the potential $V(q,t)$ identified as
\begin{equation}
      V(q,t) = \frac{1}{2} q^2 - 2  \left( \sum_{j=1}^d \Gamma_0\left(\lambda_j\right)^2 \right) \log( \cosh( e^{-t} q ) ) \;.
\end{equation}
Since the backward process is initiated around $x=0$, the speciation happens at the time $t_S$ for which the curvature of the potential changes at $x=0$. Solving $\pdv[2]{q} V(0,t_S) = 0$, we obtain:
\begin{equation}
    \label{eq:t_spec}
    t_S = \frac{1}{2}\log \left(2\sum_{j=1}^d \Gamma_0\left(\lambda_j\right)^2\right) \;.
\end{equation}
This result generalizes the one by~\cite{biroli_dynamical_2024}. The effects of the non-linearity and the manifold are encapsulated in the function $\sum_{j=1}^d \Gamma_0\left(\lambda_j\right)^2$. We proceed in order to extract the dominant behavior of this function. 

Defining the matrix  $M = \begin{bmatrix}\mu\;\vert\;\mu\;\vert\;\dots\;\vert\;\mu\end{bmatrix}\in\R^{p\times d}$, where $\mu\in\R^p$ is repeated $d$ times as columns, the sum over the functions $\Gamma_0$ can be rewritten as:
\begin{equation}
    \label{eq:gamma0tr}
    \sum_{j=1}^d \Gamma_0\left(\lambda_j\right)^2 = \frac{1}{d} \tr[\Gamma_0\left(\frac{FM}{\sqrt{p}}\right) \Gamma_0\left(\frac{FM}{\sqrt{p}}\right)^\top ] \;.
\end{equation}

The matrix $F$ is assumed to be a random matrix with i.i.d standard Gaussian entries. We make use of the Gaussian Equivalence Principle~\cite{gerace_generalisation_2020,goldt_gaussian_2022,hu_universality_2023} and write the following equivalence for $\Gamma_0\left(\frac{FM}{\sqrt{p}}\right)$:
\begin{equation}
    \label{eq:gep}
    U = \varrho_0 \one_d \one_d^\top +  \varrho_1 \left(\nicefrac{F}{\sqrt{p}} \right)   M 
    +  \varrho_* \Xi \;,
\end{equation}
where $\one_d$ is the all ones vector in $\R^d$, $\Xi\in\R^{d\times d}$ a random matrix with entries $\stackrel{i.i.d.}{\sim}\Gauss(0,1)$ and 
\begin{subequations}
\begin{align}
 \varrho_0 &=\E_{u\sim\Gauss(0,1)}\left[\Gamma_0 (u)\right] \;, \varrho_1 =\E_{u\sim\Gauss(0,1)}\left[\Gamma_0 (u)u\right] \;, \\
 \varrho_*^2  &=\E_{u\sim\Gauss(0,1)}\left[\Gamma_0 (u)^2\right] - \varrho_0^2 - \varrho_1^2  \;.
\end{align}
\end{subequations}
Since $\Gamma_0$ is an odd function, $\varrho_0 = 0$. Eventually using standard properties of Wishart matrices $F^TF/p$,  from Eq.~\eqref{eq:gep}, we obtain for $p$, $d$ large:
\begin{equation}
    \sum_{j=1}^d \Gamma_0\left(\lambda_j\right)^2 \to  \left(\varrho_1^2 \; p\right)  \norm{\tilde{\mu}}^2 + \varrho_*^2 \;,
\end{equation} 
where we have defined the rescaled mean $\tilde{\mu}=\mu/\sqrt{p}$ such that $\tilde{\mu}_j^2 \sim \bigO 1$ for $j=1, \dots, p$.
The speciation time is then:
\begin{equation}\label{eq:TS}
    t_S \approx \frac{1}{2} \log\left[ 2 \left(  \varrho_1^2  \beta d  \norm{\tilde{\mu}}^2  + \varrho_*^2 \right)  \right] \;.
\end{equation}

Therefore, the effect of an odd non-linearity on the scaling of the speciation time is a multiplicative factor given by $\varrho_1^2$, which is a positive finite number. 

\subsection{Data in a hyperplane}

If one considers a linear manifold, $\phi(y)=y$, it is straightforward to verify that $\varrho_1 = 1$ and $\varrho_* =0$. The result is then analogous to the one obtained in~\cite{biroli_dynamical_2024}, though it scales with the log of the hidden dimension $p$ instead of the dimension $d$ from the observed data. If additionally, there is no manifold, $d=p$ and $F / \sqrt{p} = I_d$, the result $t_S = (\nicefrac{1}{2}) \text{log}(2\norm{\tilde{\mu}}^2 d)$ of~\cite{biroli_dynamical_2024} is recovered.

\section{Collapse time}\label{sec:Ctime}

Assume that we run the backward process in \eqref{backward} using the empirical score function $s^e$ instead of the actual score $s$. In this case, the backward process will have probability distribution $P_t^e$. Because $a_t \to 1$ and $h_t\to 0$ as $t\to 0$, it will collapse to one of the training samples at $t=0$. Hence, we expect that, as time decreases, there exists a collapse time $t_C$ at which the trajectory is attracted to one of the training samples. To compute $t_C$, we use an analogy with the Random Energy Model (REM) of spin glass theory valid in the regime of $n=e^{\alpha d}$ training samples, first introduced in the context of diffusion models in~\cite{biroli_dynamical_2024}. Here we proceed similarly, but due to the non-linearity of the manifold model, an analogy is made with the free energy of generalized linear models (GLMs)~\cite{barbier_optimal_2019} as well.

\subsection{Reduction to a Bayesian optimal inference problem}

We consider an arbitrary sample $x_1= \phi(\frac{F\xi_1}{\sqrt p})$ where $\xi_1$ is generated from the Gaussian $q_+$, and study the distribution $P_t^e$ around this point. Let $x=a_t\phi(\frac{F\xi_1}{\sqrt p})+\sqrt{h_t}z$ be the point obtained by running the forward diffusion till a {\it small} time $t$ starting at $x_1$ (so $x$ is close to $x_1$). 
We have
\begin{align}
    P_t^e(x)  & = n^{-1}(2\pi h_t)^{-d/2}\left(e^{-\frac{\norm{z}^2}{2}}+\sum_{i=2}^ne^{-\frac{\norm{x-a_t\phi\left(\frac{F\xi_i}{\sqrt p}\right)}^2}{2h_t}}\right)
    \nonumber \\ &
    =
    n^{-1}(2\pi h_t)^{-d/2}(\Zpart_1 (t) + \Zpart_2(t)) \;.
\end{align}
We want to find the time $t_C$ such that when $t<t_C$, $\Zpart_1$ dominates over $\Zpart_2$. In other words, the score function acts as a potential well $\Vert x - a_t x_1\Vert^2/2h_t$ in which the backward trajectory "falls" towards $x_1$. 
Note that $\Zpart_1 \approx e^{-d/2}$. For the second term we have $\Zpart_2(t) = \Zpart_2^+(t) + \Zpart_2^-(t)$, where $\Zpart_2^{\pm}$ correspond to the samples generated from the Gaussian $q_{\pm}$ (there are roughly $n/2$ samples for each term).
Defining $\fe_2^\pm (t) = \lim_{d\to +\infty}\frac{1}{d}\log \Zpart_2^\pm(t)$, we expect $\Zpart_2(t) = \frac{1}{2}e^{d \fe_2^+(t)} + \frac{1}{2} e^{d \fe_2^-(t)}$ for large $d$. We shall argue that $\fe_2^+(t) > \fe_2^-(t)$, and therefore $\Zpart_2(t) \approx \frac{1}{2}e^{d \fe_2^+(t)}$ for large $d$ (this asymmetry arises because $x$ is close to $x_1$ generated from $q_+$). The collapse time can then be found from 
$e^{-d/2} \approx \frac{1}{2}e^{d \fe_2^+(t_C)}$ for large $d$, which gives the condition 
$\fe_2^+(t_C) = - \frac{1}{2}$. 

We expect to have the concentration property $\fe_2^\pm(t_C) = \lim_{d\to +\infty}\frac{1}{d}\E_x[\log \Zpart_2^\pm(t)]$ where the expectation is over $x=a_t\phi(\frac{F\xi_1}{\sqrt p})+\sqrt{h_t}z$ with probability distribution $P_t^+$, where we define
\begin{align}\label{eqn:prob_dist_pm}
    P_t^\pm(x) \equiv {(2\pi h_t)^{-d/2}}\int_{\R^p} \dd\xi\; q_\pm(\xi)e^{-\frac{\norm{x-a_t\phi(F\xi/\sqrt D)}^2}{2h_t}}. 
\end{align}
Approximating $\Zpart_2^+\approx n(2\pi h_t)^{d/2}P_t^+(x)$ (see Appendix for the validity of this approximation) in the regime
$n=e^{\alpha d}$ for large $d$, we see that
$t_C$ can be computed as the solution of the following equation:
\begin{align}\label{coll}
\alpha + \frac{1}{2}\log(2\pi h_{t_C}) + \lim_{d\to \infty}\frac{1}{d} \shortexpect_{x\sim P_{t_C}^+}\log P_{t_C}^+(x) = - \frac{1}{2} \;.
\end{align}

To justify $\fe_2^+(t)> \fe_2^-(t)$, we proceed as above and recognize that this inequality boils down to
$\shortexpect_{x\sim P_t^+}\log P_t^+(x) > \shortexpect_{x\sim P_t^+} \log P_t^-(x).$
This is indeed true because of the positivity of the Kullback-Leibler divergence between distributions $P_t^\pm$.

Now we analyze Eq.~\eqref{coll}. For non-linear $\phi$, we assume that $F$ has i.i.d. $\cN{0,1}$ entries. First, we notice that without loss of generality, we can assume that $\mu_+ = m\one_p$, where $\one_p$ is the all ones vector of dimension $p$ and $m=\norm{\mu_+}/\sqrt{p}$. This can be seen by rotating the axis of integration in \eqref{eqn:prob_dist_pm}.  We then recognize on the left hand side the Bayesian optimal free energy of a GLM. This is an inference model where we have observations
$x = a_t\phi(F\xi_1/\sqrt{p})+\sqrt{h_t} z$, with $\xi_1$ a signal to be estimated, $z$ Gaussian additive noise, and $a_t^2/h_t= e^{-2t}/(1-e^{-2t})$ the signal-to-noise ratio. When the Bayesian statistician uses the "correct" prior probability distribution $q_+$, the log-normalizing factor of the posterior distribution is precisely the free energy on the left hand side of \eqref{coll}. This is a statistical mechanics spin-glass problem with Nishimori symmetry, whose 
rigorous theory was developed in~\cite{barbier_optimal_2019}. Note that the other free energy (corresponding to $\Zpart_2^-$) does not satisfy Nishimori symmetry because it corresponds to a mismatched prior $q_-$ used by the statistician. 
We have:
\begin{align*}
    \lim_{p\to \infty} \frac{1}{p} \shortexpect_{x\sim P_t^+}\log P_t^+(x) &= \sup_{q\in[0,\rho+m^2]} \inf_{r\ge 0} f_{\text{RS}}(q,r)\eqdef f^\star(t)
    \label{eq:fstar}
\end{align*}
where
\begin{equation}\label{FRS}
    f_{\text{RS}}(q,r) = \psi(r) + \beta^{-1}\Psi(q)-rq/2 \;,
\end{equation}
\begin{equation*}
    \psi(r) = \shortexpect_{X_0,Z_0}{\log \int \dd w\frac{e^{-\frac{(w-m)^2}{2\rho}}}{\sqrt{2\pi\rho}}e^{rwX_0+\sqrt{r}xZ_0-rx^2/2}} \;,
\end{equation*}
\begin{equation*}
    \Psi(q) = \shortexpect_{Y_0,V}{\log \int \dd w \frac{e^{-w^2/2}}{\sqrt{2\pi}}\frac{e^{\frac{(Y_0-a_t\phi(\sqrt{q}V+\sqrt{m^2+\rho-q}w))^2}{2h_t}}}{\sqrt{2\pi h_t}}} \;,
\end{equation*}
with $X_0\sim\cN{m,\rho}$, $Z_0,V,W,Z\sim\cN{0,1}$ and $Y_0 = a_t\phi(\sqrt{q}V+\sqrt{m^2+\rho-q}W)+\sqrt{h_t}Z$. It is direct to compute $\psi$. We get $\psi(r) = \frac{r(m^2+\rho)}{2}-\frac{1}{2}\log(1+r\rho)$.

\subsection{General manifold Model}
\begin{figure}[htbp]
    \centering
    \includegraphics[scale=0.5]{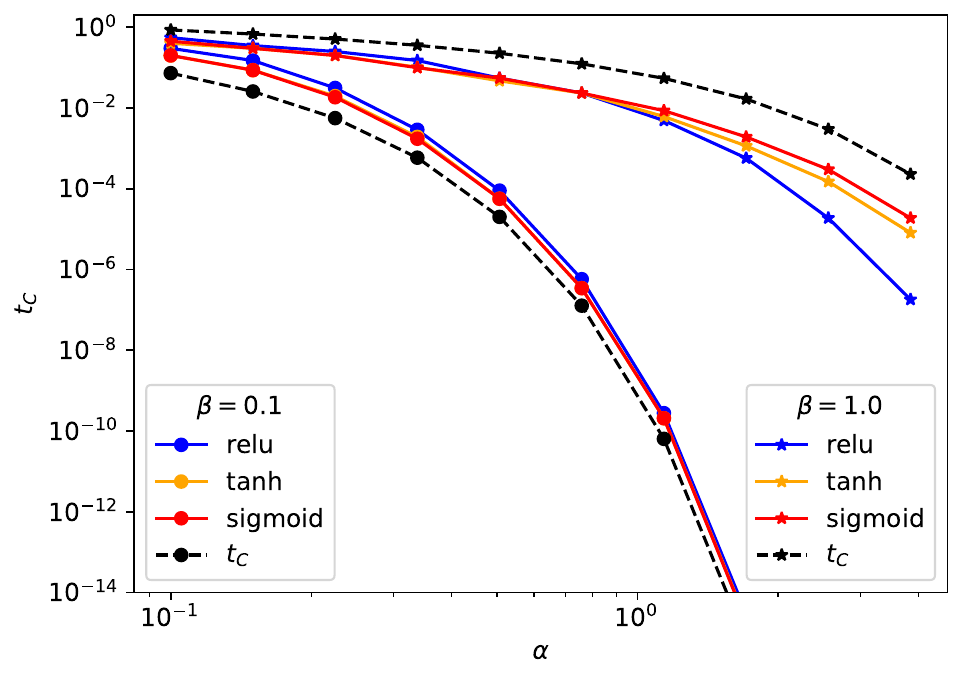}
    \caption{Collapse time for different non-linearities. The curve $t_C$ refers to the collapse time obtained using (\ref{eqn:tc_linear_manifold}).}
    \label{fig:tc_nonlin_manifold}
\end{figure}

For a nonlinear activation $\phi$, $\Psi$ needs to be computed numerically. Finally, the collapse time is found by solving 
\begin{align}\label{eq:TC}
\alpha + (\nicefrac{1}{2})\log(2\pi h_{t_C}) + \beta f^\star(t_C)=-\nicefrac{1}{2}.
\end{align}
Fig.~\ref{fig:tc_nonlin_manifold} illustrates the collapse time obtained for \textit{relu}, \textit{tanh} and \textit{sigmoid} non-linearities. 

\subsection{Data in a hyperplane}
\begin{figure}
    \centering
    \includegraphics[scale=0.5]{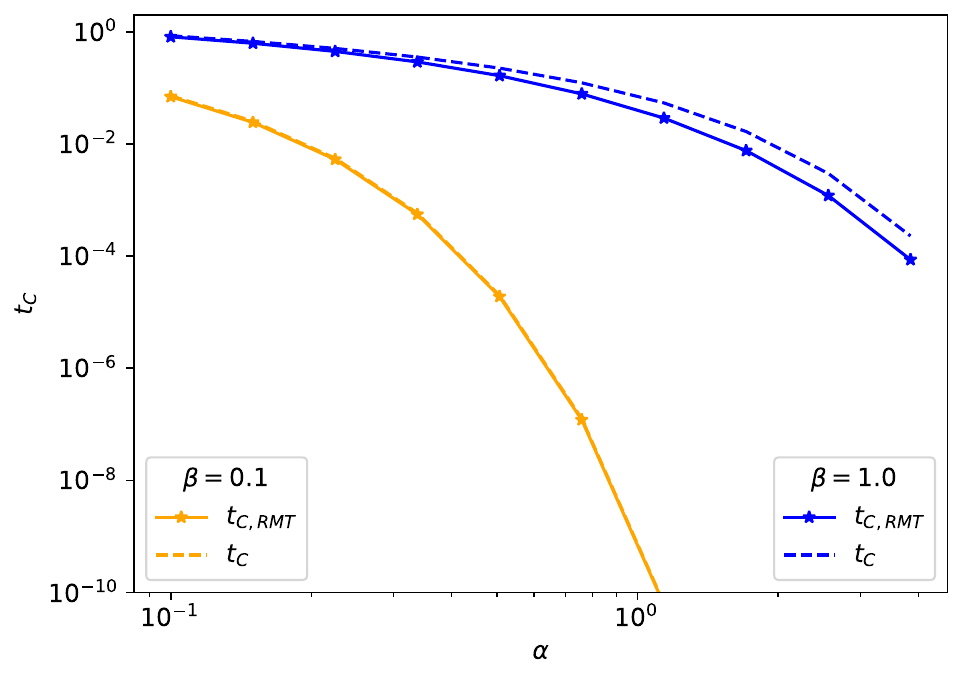}
    \caption{Collapse time for linear manifold. $t_{c,RMT}$ for random $F$ and $t_C$ for isometric $F$.  (\ref{eqn:tc_linear_manifold})}
    \label{fig:tc_linear_rmt}
\end{figure}

For a linear activation function $\phi(u)=u$, the data lies in a hyperplane of dimension $p<d$. In this case we can compute the Gaussian integral which yields (we set $\eta_t \equiv a_t^2/h_t)$), $P_t^+(x) = ( (2\pi)^d \det\Sigma)^{-1/2} e^{-\frac{1}{2}(x-\mu_+)^\top\Sigma^{-1}(x-\mu_+)} $. Therefore, $\fe_2^+$ is given by
\begin{align*}
    \fe_2^+(t) = \alpha + \frac{1}{2}\log(h_t)  - \lim_{d\to +\infty} \frac{1}{2d}\log \det\Sigma-\frac{1}{2} \;,
\end{align*}
where $\Sigma = h_t(\eta_t FF^T+I_d)$. For the determinant, we have $\frac{1}{d}\log\det\Sigma =  \log(h_t) + \frac{1}{d}\log\det(\eta_t\frac{FF^T}{p}+I_d)$. Thus, we find collapse time by the condition 
\begin{equation}\label{eq:collapsetimecondition_linear}
    \alpha-\lim_{d\to +\infty}\frac{1}{2d}\log\det(\eta_{t_C}FF^T/p+I_d)=0.
\end{equation}
Now, we specialize to the cases of random and deterministic isometric matrices for $F$. The collapse times are compared on Fig.~\ref{fig:tc_linear_rmt}. As expected, the differences are negligible for small $\beta$. 

\subsubsection{Deterministic isometry for $F$}
Because $F^TF/p = I_p$, the $d\times d$ matrix $FF^T/p$ has $p$ eigenvalues equal to $1$ and $d-p$ eigenvalues equal to $0$. With this remark we can compute $\frac{1}{d}\log\det(\eta_tFF^T/p+I_d) = \beta\log(1+\eta_t)$, which yields 
\begin{equation}\label{eqn:tc_linear_manifold}
    t_C = (\nicefrac{1}{2})\log\left(1+(e^{2\alpha/\beta}-1)^{-1}\right).
\end{equation}
When $d=p$, that is $\beta=1$, we recover the formula from \cite{biroli_dynamical_2024}.

\subsubsection{Random matrix for $F$}
Using the Marchenko-Pastur distribution and standard techniques for the random matrix $FF^T/p$ in the large dimensional limit, we obtain $\frac{1}{d}\log\text{det}(\eta_tFF^T/p+I_d) = \beta\log\left(1+\frac{\eta_t}{\beta}-\frac{1}{4}h(\frac{\eta_t}{\beta},\beta)\right)+\log\left(1+\eta_t-\frac{1}{4}h(\frac{\eta_t}{\beta},\beta)\right)-\frac{\beta}{4\eta_t}h(\frac{\eta_t}{\beta},\beta)$, where $h(x,z) = \left(\sqrt{x(1+\sqrt{z})+1}-\sqrt{x(1-\sqrt{z})+1}\right)^2$.
The condition \eqref{eq:collapsetimecondition_linear} to find the collapse time is now a closed equation which can be easily solved numerically. 

This case can also be treated through the theory for a general manifold by computing explicitly
$\Psi(q) = -\frac{1}{2}-\frac{1}{2}\log(2\pi(a_t^2(m^2+\rho-q)+h_t))$. From there we deduce $\fe_2^+(t)$. The result agrees with the random matrix calculation.
 
\section{Conclusion}
From our results on speciation and collapse time, we conclude that these times are much smaller when the data comes from a low-dimensional manifold. In particular, the number of samples required to keep these times at $\bigO{1}$ scales as $\bigO{e^p}$ for manifold data, where $p$ is the dimension of the manifold. This is advantageous, as we need these times to be as small as possible to mitigate memorization.
Obviously, it would be desirable to generalize the analysis to more general data models on manifolds. Even for Gaussian mixtures with more than two centers the situation can become complicated, with potentially many speciation and collapse times depending on the location of the centers. Furthermore, it would be desirable establish the analysis presented here on mathematically rigorous grounds.

\appendix
In Eq.~\eqref{coll}, we used the approximation $\fe_2^+\approx \alpha+ \frac{1}{2}\log(2\pi h_t) +\lim_{d\to\infty}\shortexpect_x\log P_t^+(x)$. This approximation is however delicate and is valid only when $t$ is large. To obtain $\fe_2^+$ for all $t$, we can view it as the log partition function of a REM \cite{mezard_information_2009, lucibello_exponential_2024, biroli_dynamical_2024}. Let 
\begin{equation*}
    P_{t,\lambda}^+(x) = (2\pi h_t)^{-d/2}\int_{\R^p} \dd\xi\; q_+(\xi) e^{-\lambda\frac{\norm{x-a_t\phi(F\xi/\sqrt p)}^2}{2h_t}},
\end{equation*}
and $g_t(\lambda) = \lim_{d\to\infty}\frac{1}{d}\shortexpect_x\log P_t^+(x,\lambda)$. Then, by REM theory, the function $\fe_2^+$ undergoes a \textit{condensation} phase transition at time $t^*$. The time $t^*$ can be obtained by the condition $\alpha_n+g_{t^*}(1)-g_{t^*}'(1)=0$. For $t \ge t^*$, $\alpha_n+g_{t}(1)$ well approximates $\fe_2^+$. This is not the case however for $t<t^*$. Nevertheless, from the following argument, we find that $g_t'(1) = -1/2$:
\begin{align*}
        &-g_t'(1) = -\lim_{d\to \infty}\frac{1}{d}\bE{x}{\frac{1}{P_t^+(x)}\frac{\partial P_{t,\lambda}(x)}{\partial \lambda}\Big|_{\lambda=1}}\\
        &
        = \lim_{d\rightarrow \infty}\frac{1}{d}\shortexpect_x \expectCond{\frac{\norm{x-a_tx_1}^2}{2h_t}}{x}
        % &= -\lim_{d\rightarrow \infty}\frac{1}{d}\expect{\frac{\norm{x-a_tx_0}^2}{2h_t}}
        = \lim_{d\rightarrow \infty}\frac{1}{2d}\shortexpect\norm{z}^2 = \frac{1}{2}.        
    \end{align*}
This implies that $t^*$ and $t_c$ calculated using \eqref{coll} are the same. Thus, the approximation made in \eqref{coll} is valid for $t\ge t_c$.

\section*{Acknowledgment}
The work of A. J. G and R. V has been supported by Swiss National Science Foundation grant number 200021-204119.

\bibliographystyle{IEEEtran}
\bibliography{references.bib}

\begin{thebibliography}{10}
\providecommand{\url}[1]{#1}
\csname url@samestyle\endcsname
\providecommand{\newblock}{\relax}
\providecommand{\bibinfo}[2]{#2}
\providecommand{\BIBentrySTDinterwordspacing}{\spaceskip=0pt\relax}
\providecommand{\BIBentryALTinterwordstretchfactor}{4}
\providecommand{\BIBentryALTinterwordspacing}{\spaceskip=\fontdimen2\font plus
\BIBentryALTinterwordstretchfactor\fontdimen3\font minus \fontdimen4\font\relax}
\providecommand{\BIBforeignlanguage}[2]{{%
\expandafter\ifx\csname l@#1\endcsname\relax
\typeout{** WARNING: IEEEtran.bst: No hyphenation pattern has been}%
\typeout{** loaded for the language `#1'. Using the pattern for}%
\typeout{** the default language instead.}%
\else
\language=\csname l@#1\endcsname
\fi
#2}}
\providecommand{\BIBdecl}{\relax}
\BIBdecl

\bibitem{sohl-dickstein_deep_2015}
J.~Sohl-Dickstein, E.~Weiss, N.~Maheswaranathan, and S.~Ganguli, ``\BIBforeignlanguage{en}{Deep {Unsupervised} {Learning} using {Nonequilibrium} {Thermodynamics}},'' in \emph{\BIBforeignlanguage{en}{Proceedings of the 32nd {International} {Conference} on {Machine} {Learning}}}.\hskip 1em plus 0.5em minus 0.4em\relax PMLR, Jun. 2015, pp. 2256--2265, iSSN: 1938-7228.

\bibitem{song_generative_2019}
Y.~Song and S.~Ermon, ``Generative {Modeling} by {Estimating} {Gradients} of the {Data} {Distribution},'' in \emph{Advances in {Neural} {Information} {Processing} {Systems}}, vol.~32.\hskip 1em plus 0.5em minus 0.4em\relax Curran Associates, Inc., 2019.

\bibitem{ho_denoising_2020}
J.~Ho, A.~Jain, and P.~Abbeel, ``Denoising {Diffusion} {Probabilistic} {Models},'' in \emph{Advances in {Neural} {Information} {Processing} {Systems}}, vol.~33.\hskip 1em plus 0.5em minus 0.4em\relax Curran Associates, Inc., 2020, pp. 6840--6851.

\bibitem{yang_diffusion_2024}
L.~Yang, Z.~Zhang, Y.~Song, S.~Hong, R.~Xu, Y.~Zhao, W.~Zhang, B.~Cui, and M.-H. Yang, ``\BIBforeignlanguage{en}{Diffusion {Models}: {A} {Comprehensive} {Survey} of {Methods} and {Applications}},'' \emph{\BIBforeignlanguage{en}{ACM Computing Surveys}}, vol.~56, no.~4, pp. 1--39, Apr. 2024.

\bibitem{anderson_reverse-time_1982}
B.~D.~O. Anderson, ``Reverse-time diffusion equation models,'' \emph{Stochastic Processes and their Applications}, vol.~12, no.~3, pp. 313--326, May 1982.

\bibitem{haussmann_time_1986}
U.~G. Haussmann and E.~Pardoux, ``Time {Reversal} of {Diffusions},'' \emph{The Annals of Probability}, vol.~14, no.~4, pp. 1188--1205, Oct. 1986, publisher: Institute of Mathematical Statistics.

\bibitem{yoon_diffusion_2023}
T.~Yoon, J.~Y. Choi, S.~Kwon, and E.~K. Ryu, ``\BIBforeignlanguage{en}{Diffusion {Probabilistic} {Models} {Generalize} when {They} {Fail} to {Memorize}},'' in \emph{\BIBforeignlanguage{en}{{ICML} 2023 {Workshop} on {Structured} {Probabilistic} {Inference} and {Generative} {Modeling}}}, Jul. 2023.

\bibitem{biroli_generative_2023}
G.~Biroli and M.~Mézard, ``Generative diffusion in very large dimensions,'' \emph{Journal of Statistical Mechanics: Theory and Experiment}, vol. 2023, no.~9, p. 093402, Sep. 2023.

\bibitem{biroli_dynamical_2024}
G.~Biroli, T.~Bonnaire, V.~de~Bortoli, and M.~Mézard, ``\BIBforeignlanguage{en}{Dynamical regimes of diffusion models},'' \emph{\BIBforeignlanguage{en}{Nature Communications}}, vol.~15, no.~1, p. 9957, Nov. 2024, publisher: Nature Publishing Group.

\bibitem{raya_spontaneous_2023}
G.~Raya and L.~Ambrogioni, ``Spontaneous symmetry breaking in generative diffusion models,'' in \emph{Advances in {Neural} {Information} {Processing} {Systems}}, vol.~36.\hskip 1em plus 0.5em minus 0.4em\relax Curran Associates, Inc., 2023, pp. 66\,377--66\,389.

\bibitem{li_critical_2024}
M.~Li and S.~Chen, ``\BIBforeignlanguage{en}{Critical windows: non-asymptotic theory for feature emergence in diffusion models},'' in \emph{\BIBforeignlanguage{en}{Proceedings of the 41st {International} {Conference} on {Machine} {Learning}}}.\hskip 1em plus 0.5em minus 0.4em\relax PMLR, Jul. 2024, pp. 27\,474--27\,498, iSSN: 2640-3498.

\bibitem{sclocchi_phase_2025}
A.~Sclocchi, A.~Favero, and M.~Wyart, ``\BIBforeignlanguage{en}{A phase transition in diffusion models reveals the hierarchical nature of data},'' \emph{\BIBforeignlanguage{en}{Proceedings of the National Academy of Sciences}}, vol. 122, no.~1, p. e2408799121, Jan. 2025.

\bibitem{goldt_modeling_2020}
S.~Goldt, M.~Mézard, F.~Krzakala, and L.~Zdeborová, ``\BIBforeignlanguage{en}{Modeling the {Influence} of {Data} {Structure} on {Learning} in {Neural} {Networks}: {The} {Hidden} {Manifold} {Model}},'' \emph{\BIBforeignlanguage{en}{Physical Review X}}, vol.~10, no.~4, p. 041044, Dec. 2020.

\bibitem{hand_global_2020}
P.~Hand and V.~Voroninski, ``Global {Guarantees} for {Enforcing} {Deep} {Generative} {Priors} by {Empirical} {Risk},'' \emph{IEEE Transactions on Information Theory}, vol.~66, no.~1, pp. 401--418, Jan. 2020, conference Name: IEEE Transactions on Information Theory.

\bibitem{luneau_tensor_2020}
C.~Luneau and N.~Macris, ``Tensor {Estimation} {With} {Structured} {Priors},'' \emph{IEEE Journal on Selected Areas in Information Theory}, vol.~1, no.~3, pp. 705--722, Nov. 2020, conference Name: IEEE Journal on Selected Areas in Information Theory.

\bibitem{achilli_losing_2024}
B.~Achilli, E.~Ventura, G.~Silvestri, B.~Pham, G.~Raya, D.~Krotov, C.~Lucibello, and L.~Ambrogioni, ``Losing dimensions: {Geometric} memorization in generative diffusion,'' Oct. 2024, arXiv:2410.08727.

\bibitem{ventura_manifolds_2024}
E.~Ventura, B.~Achilli, G.~Silvestri, C.~Lucibello, and L.~Ambrogioni, ``Manifolds, {Random} {Matrices} and {Spectral} {Gaps}: {The} geometric phases of generative diffusion,'' Oct. 2024, arXiv:2410.05898.

\bibitem{gardiner_stochastic_2009}
C.~W. Gardiner, \emph{\BIBforeignlanguage{eng}{Stochastic methods: a handbook for the natural and social sciences}}, 4th~ed., ser. Springer series in synergetics.\hskip 1em plus 0.5em minus 0.4em\relax Berlin Heidelberg: Springer, 2009, no.~13.

\bibitem{song_maximum_2021}
Y.~Song, C.~Durkan, I.~Murray, and S.~Ermon, ``Maximum {Likelihood} {Training} of {Score}-{Based} {Diffusion} {Models},'' in \emph{Advances in {Neural} {Information} {Processing} {Systems}}, vol.~34.\hskip 1em plus 0.5em minus 0.4em\relax Curran Associates, Inc., 2021, pp. 1415--1428.

\bibitem{barbier_optimal_2019}
J.~Barbier, F.~Krzakala, N.~Macris, L.~Miolane, and L.~Zdeborová, ``\BIBforeignlanguage{en}{Optimal errors and phase transitions in high-dimensional generalized linear models},'' \emph{\BIBforeignlanguage{en}{Proceedings of the National Academy of Sciences}}, vol. 116, no.~12, pp. 5451--5460, Mar. 2019.

\bibitem{kibble_extension_1945}
W.~F. Kibble, ``\BIBforeignlanguage{en}{An extension of a theorem of {Mehler}'s on {Hermite} polynomials},'' \emph{\BIBforeignlanguage{en}{Mathematical Proceedings of the Cambridge Philosophical Society}}, vol.~41, no.~1, pp. 12--15, Jun. 1945.

\bibitem{gerace_generalisation_2020}
F.~Gerace, B.~Loureiro, F.~Krzakala, M.~Mezard, and L.~Zdeborova, ``\BIBforeignlanguage{en}{Generalisation error in learning with random features and the hidden manifold model},'' in \emph{\BIBforeignlanguage{en}{Proceedings of the 37th {International} {Conference} on {Machine} {Learning}}}.\hskip 1em plus 0.5em minus 0.4em\relax PMLR, Nov. 2020, pp. 3452--3462, iSSN: 2640-3498.

\bibitem{goldt_gaussian_2022}
S.~Goldt, B.~Loureiro, G.~Reeves, F.~Krzakala, M.~Mezard, and L.~Zdeborova, ``\BIBforeignlanguage{en}{The {Gaussian} equivalence of generative models for learning with shallow neural networks},'' in \emph{\BIBforeignlanguage{en}{Proceedings of the 2nd {Mathematical} and {Scientific} {Machine} {Learning} {Conference}}}.\hskip 1em plus 0.5em minus 0.4em\relax PMLR, Apr. 2022, pp. 426--471, iSSN: 2640-3498.

\bibitem{hu_universality_2023}
H.~Hu and Y.~M. Lu, ``Universality {Laws} for {High}-{Dimensional} {Learning} {With} {Random} {Features},'' \emph{IEEE Transactions on Information Theory}, vol.~69, no.~3, pp. 1932--1964, Mar. 2023, conference Name: IEEE Transactions on Information Theory.

\bibitem{mezard_information_2009}
M.~Mézard and A.~Montanari, Eds., \emph{Information, {Physics}, and {Computation}}.\hskip 1em plus 0.5em minus 0.4em\relax Oxford University Press, Jan. 2009.

\bibitem{lucibello_exponential_2024}
C.~Lucibello and M.~Mézard, ``\BIBforeignlanguage{en}{The {Exponential} {Capacity} of {Dense} {Associative} {Memories}},'' \emph{\BIBforeignlanguage{en}{Physical Review Letters}}, vol. 132, no.~7, p. 077301, Feb. 2024, arXiv:2304.14964 [cond-mat].

\end{thebibliography}

\end{document}